\documentclass[11pt]{article}
\usepackage{latexsym}
\usepackage[dvips]{epsfig}
\newcommand{\qed}{\mbox{$\Diamond$}\vspace{\baselineskip}}

\newtheorem{theorem}{Theorem}[section]
\newtheorem{proposition}[theorem]{Proposition}
\newtheorem{lemma}[theorem]{Lemma}
\newtheorem{definition}[theorem]{Definition}
\newtheorem{corollary}[theorem]{Corollary}
\newtheorem{example}[theorem]{Example}

\newtheorem{conjecture}[theorem]{Conjecture} 
\newenvironment{proof}{\noindent {\bf Proof:}}{{\qed}}

\newcommand{\vanish}[1]{} 
\begin{document}
\title{The limit of a Stanley-Wilf sequence is
not always rational, and layered patterns beat monotone patterns}

\author{Mikl\'os B\'ona
\thanks{University of Florida, Gainesville FL 32611-8105. Partially supported
by an NSA Young Investigator Award.
 Email: bona@math.ufl.edu.}
}

 \date{}

\maketitle

\begin{abstract} We show the first known example for a pattern $q$ for which 
$L(q)=\lim_{n\rightarrow \infty} \sqrt[n]{S_n(q)}$ is not an integer. We 
find the exact value of the limit and show that it is irrational, 
but algebraic. Then we
generalize our results to an infinite sequence of patterns. We
provide further generalizations that start explaining why certain patterns
are easier to avoid than others. Finally, we show that if $q$ is a 
layered pattern of length $k$, then $L(q)\geq (k-1)^2$ holds. 
\end{abstract}

\section{Introduction}
 Let $S_n(q)$ be the number of permutations of length $n$ 
(or, in what follows,
$n$-permutations) that avoid the pattern $q$. For a brief introduction to
the area of pattern avoidance, see \cite{bonawalk}; for a more detailed
introduction, see \cite{Bona}. A recent spectacular result of Marcus
and Tardos \cite{marcus} shows that for any pattern $q$, there exists
a constant $c_q$ so that $S_n(q)<c_q^n$ holds for all $n$. As pointed out
by Arratia in \cite{arratia}, this is equivalent to the statement 
that $L(q)=\lim_{n\rightarrow \infty} \sqrt[n]{S_n(q)}$ exists. Let us call
the sequence $ \sqrt[n]{S_n(q)}$ a {\em Stanley-Wilf} sequence. 
It is a natural 
and intriguing question to ask what the limit $L(q)$
 of a Stanley-Wilf sequence
 can be, for various patterns
$q$. 

The main reason this question has been so intriguing is that in all cases
where  $L(q)$ has been known, it has
been known to be an {\em integer}. Indeed, the results previously known
are listed below. 

\begin{enumerate}
\item When $q$ is of length three, then 
 $L(q)=4$. This follows from the
well-known fact \cite{simion} that in this case, $S_n(q)={2n\choose n}/(n+1)$.
\item When $q=123\cdots k$, or when $q$ is such that $S_n(q)=S_n(12\cdots k)$,
then  $L(q)=(k-1)^2$. This follows
from an asymptotic formula of Regev \cite{regev}.
\item When $q=1342$, or when  $q$ is such that $S_n(q)=S_n(1342)$, then
 $L(q)=8$. See \cite{bona} for this
result and an exact formula for the numbers $S_n(1342)$.
\end{enumerate}

In this paper, we show that $L(q)$
is {\em not} always an integer. We achieve this by proving that 
$14<L(12453)<15$. Then we compute
the exact value of this limit, and see that it is not even rational; it is
the number $9+4\sqrt{2}$. We compute the limit of the Stanley-Wilf
sequence for an infinite sequence of patterns, and see that as the length
$k$ of these patterns grows,  $L(q)$
will fall further and further below the largest known possible value,
$(k-1)^2$. Finally, we show that while for certain patterns, our methods 
provide the exact value of the limit of the Stanley-Wilf sequence, for 
certain others they only provide a {\em lower bound} on this limit. This
starts explaining why certain patterns are easier to avoid than others.
Among other results, we will confirm a seven-year old conjecture  by proving
that in the sense of logarithmic asymptotics,
a layered pattern $q$ is always easier to avoid than the monotone pattern of
the same length. 

\section{Proving an upper bound}
Let $p=p_1p_2\cdots p_n$ be a permutation. Recall that $p_i$ is called a
{\em left-to-right minimum} of $p$ if $p_j>p_i$ for all $j<i$. In other words,
a left-to-right minimum is an entry that is smaller than everything on its 
left. Note that $p_1$ is always a left-to-right minimum, and so is the entry
1 of $p$. Also note that the left-to-right minima of $p$ always form a
decreasing sequence. For the rest of this paper, entries that are not
left-to-right minima are called {\em remaining entries}.

Now we are in a position to prove our promised upper bound for the
numbers $S_n(12453)$.

\begin{lemma} \label{lupper} For all positive integers $n$, we have 
\[S_n(12453)< (9+4\sqrt{2})^n<14.66^n.\]
\end{lemma}

\begin{proof}
Let $p$ be a permutation counted by $S_n(12453)$, and let $p$ have $k$
left-to-right minima. Then we have at most ${n\choose k}$ choices for the
set of these left-to-right minima, and we have at most ${n\choose k}$
choices for their positions. The string of the remaining entries has
to form a 1342-avoiding permutation of length $n-k$. Indeed, if there was a 
copy $acdb$ of 1342 among the entries that are not left-to-right minima, then
we could complete it to a 12453 pattern by simply prepending it by the closest
left-to-right minimum that is on the left of $a$. The number 
of 1342-avoiding permutations on $n-k$ elements
 is less than $8^{n-k}$ as we know from \cite{bona}.
This shows that
\begin{eqnarray*} S_n(12453) & < & \sum_{k=1}^n {n\choose k}^2\cdot 8^{n-k} \\
& < & \sum_{k=1}^n \left ({n\choose k}\cdot \sqrt{8}^{n-k}\right )^2 \leq
\left (\sum_{k=1}^n {n\choose k}\cdot \sqrt{8}^{n-k} \right )^2 \\
& < & (1+ \sqrt{8})^{2n}=(9+4\sqrt{2})^n, 
\end{eqnarray*}
and the proof is complete.
\end{proof}

\begin{corollary} \label{upper}
 We have
\[L(12453) < 14.66. \]
\end{corollary}

\section{Proving a lower bound}
We have seen in Corollary \ref{upper} 
that  $L(12453)  = 9+4\sqrt{2}<14.66$. In order to prove
 that this limit is not an integer, it suffices to show that it is larger than
14. In what follows, we are going to work towards a good lower bound for
the numbers $S_n(12453)$, and thus the number $L(12453)$.

Where is the waste in the proof of the upper bound in the previous section?
The waste is that
there are some choices for the left-to-right minima that are incompatible
with some choices for the 1342-avoiding permutation of the remaining entries.
This is a crucial concept of the upcoming proof, so we will make it more
precise. 

We have mentioned in the previous section, that determining
 the left-to-right minima of a permutation $p$ means
to determine the set $T$ of positions these minima will be, and to determine
 the set
$Z$ of entries that are the left-to-right minima. In other words, the
ordered pair $(T,Z)$ of equal-sized subsets of $[n]=\{1,2,\cdots ,n\}$
 describes the
left-to-right minima of $p$.  

\begin{definition}
Let $n$ be a positive integer, and let $m\leq n$ be a positive integer.
Let $T$ and $Z$ be two $m$-element subsets of $[n]$. Finally, let $S$ be
a permutation of the elements of the set $[n]-Z$. If there exists
an $n$-permutation $p$ so that its left-to-right minima are precisely
the elements of $Z$, they are located in positions belonging to $T$, and
its string of remaining entries is $S$, then we say that the triple
$(T,Z,S)$ is {\em compatible}. Otherwise, we say that the triple
$(T,Z,S)$ is {\em incompatible}.
\end{definition}

Clearly, if $(T,Z,S)$ is {\em compatible}, then there is exactly one 
permutation $p$ satisfying all criteria specified by $(T,Z,S)$.

\begin{example} If $n=4$, and $T=\{1,3\}$, $Z=\{1,2\}$, and $S=43$, 
then $(T,Z,S)$ is {\em compatible} as shown by the permutation 2413.
\end{example}

\begin{example}  If $n=4$, and $T=\{1,3\}$,  $Z=\{1,3\}$, and $S=24$, 
then  $(T,Z,S)$ is {\em incompatible}. Indeed, the only permutation
allowed by $T$ and $S$ is 3214, but for this permutation $Z=\{1,2,3\}$, 
not $\{1,3\}$.
\end{example}

Returning to the method by which we proved our upper bound for $L(12453)$,
 we will show that in a sufficient number of  cases, our
triples $(T,Z,S)$ are compatible. This will show that the upper bound
is quite close to the precise value of $L(12453)$.

What is a good way to check that a particular choice $(T,Z)$ of left-to-right 
minima is compatible with a particular choice of $S$? For shortness,
let us call the procedure of putting together $S$ and a string  $(T,Z)$ of
 left-to-right minima {\em merging}.
 One has to check that in the permutation obtained
by merging our left-to-right minima with $S$, the left-to-right minima
are indeed the entries in $Z$. That is, there are no additional 
left-to-right minima, and there the entries in $Z$ are indeed all 
left-to-right minima. 
 This is achieved exactly when any remaining
entry is larger than the closest left-to-right minimum on its left.

In our efforts to find a good lower bound on $L(12453)$, we will only  
 consider a special kind of permutations. 
Let $N$ be a positive integer so that $S_n(1342)>7.99^n$ for all $n>N$. 
(We know from \cite{bona} that such an $N$ exists as 
$L(1342)= 8$.)

\begin{figure}[ht]
 \begin{center}
  \epsfig{file=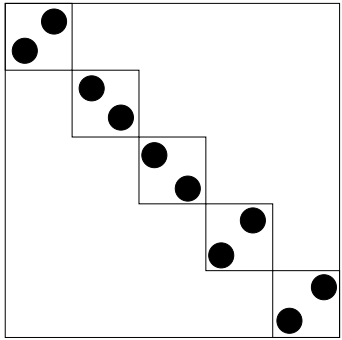}
  \caption{A block-structured string. }
  \label{blocks}
 \end{center}
\end{figure}

Consider permutations whose string $S$
 of remaining entries  has the following property. If we cut
$S$ into $\lfloor |S|/N \rfloor $ blocks of consecutive entries of length 
$N$ each (the last 
block can be of size between $N$ and $2N-1$), then the entries of any given
 block $B$ are all smaller
than the entries on any block on the left of $B$, and larger than the entries
of any block on the right of $B$. Let us call these strings $S$ 
{\em block-structured}. See Fig. \ref{blocks} for the generic diagram of
a block-structured string in the (unrealistic) case of $N=2$.

The number of such strings $S$ is obviously more
than $7.99^{\lfloor |S|/N \rfloor }$.
 It is obvious that they are all 1342-avoiding as a 
1342-pattern cannot start in a block and end in another one. 
We claim that a sufficient number of these strings
 $S$ will be compatible with a sufficient number of the choices $(T,Z)$ of
 left-to-right minima.

First, look at the very special case when $S$ is {\em decreasing}. In 
this case, we will write $S^{dec}$ instead of $S$.
Now our permutation $p$ consists of two decreasing sequences (so
it is 123-avoiding), namely  the left-to-right minima
and $S^{dec}$. The following Proposition
is very well-known. 

\begin{proposition} \label{pnara}
Let $1\leq m\leq n$.
Then the number of 123-avoiding $n$-permutations having exactly
$m$ left-to-right minima is
\begin{equation} \label{nara} A(n,m)=\frac{1}{n}{n\choose m}{n\choose m+1},
\end{equation} a {\em Narayana} number. 
\end{proposition}

For a proof, see \cite{stanley} or \cite{Bona}.

The significance of this result for us is the following. If we just wanted to
merge $(T,Z)$ and $S$ together, with no regard to the existing constraints, 
the total number of ways to do that would be of course at most ${n\choose m} 
\cdot {n\choose m}$. The above formula shows that roughly $\frac{1}{n}$
 of these
mergings will actually be good, that is, they will not violate any 
constraints, they will lead to compatible triples $(T,Z,S)$.
 The factor $\frac{1}{n}$ is not a significant loss from our point
of view, since $\lim_{n\rightarrow \infty}  \sqrt[n]{1/n}=1$.

Now let us return to the general case of block-structured  strings $S$.
In other words, take a 123-avoiding $n$-permutation $(T,Z,S^{dec})$,
 and replace its 
string $S^{dec}$ by a block-structured string $S$ {\em taken on the entries
that belong to} $S^{dec}$. We claim that after this replacement,
 a sufficient number of triples
$(T,Z,S)$ will be compatible.

Here is the outline of the proof of that claim. Because of
the structure of a block-structured $S$, it is true that every entry in $S$
 is at most
$N$ positions away from the position it was
 $S^{dec}$.
Therefore, if we merge $(T,Z)$ and $S^{dec}$ together so that each remaining 
entry $x$ is not simply larger than the  left-to-right minimum that is
 closest to
and preceding the position $j$ of $x$,
 but also larger than the left-to-right minimum
closest to and preceding 
 position $j-N$, then we will be done. Indeed, in this case
replacing $S^{dec}$ by any block-structured string $S$ will not violate
any constraints. 

 Therefore, we will have a lower bound for the number
of compatible triples $(T,Z,S)$ if we find a lower bound for the
number of  compatible triples $(T,Z,S^{dec})$  in which 
each remaining entry has the mentioned stronger property. 

In order to find such a lower bound, take a 123-avoiding permutation $p'$
which is of length $n-N$. Let $p'$ have $m$ left-to-right minima.
 Denote $(T',Z')$ the string of the left-to-right 
minima of $p'$, and let $S^{dec'}$ denote the decreasing string of remaining
entries of $p'$. Now prepend $p'$ with the decreasing string 
taken on the $N$-element 
set $\{n-N+1,n-N+2,\cdots ,n\}$, to get an $n$-permutation.
In this $n$-permutation, move each of the original $m$ left-to-right minima
of $p'$ to the left by $N$ positions.
  Let us call the obtained $n$-permutation
$p''$.

 It is then clear that the left-to-right minima of $p''$ are the same
as the left-to-right minima of $p'$. Furthermore, because of the translation
we used to create our new permutation,  $p''$ 
has the property that if $x$ is a remaining entry of $p''$ and is in 
position $j$, then $x$ is larger than the left-to-right minimum that is
closest to, and preceding, position $j-N$. 

Now we can use the argument that we outlined four paragraphs ago. For easy
reference, we sketch that argument again.
 If $S^{dec}$ is replaced by
any block-structured permutation of the same size taken on the same
set of elements, (resulting in the
$n$-permutation $p*$) then each remaining
entry $x$ will move within its block only, that is, $x$ will move at
most $N$ positions from its original position. Therefore, $x$ will still
be larger than the left-to-right minimum closest to it and preceding it.

This shows that if $p'$ and $(T',Z')$ lead to a compatible triple,
 then so too will $p*$ and $(T,Z)$, where $(T,Z)$ describes the left-to-right
minima of $p*$. Proposition \ref{pnara}
implies that the number of compatible triples 
 $(T',Z',p')$ is
$\frac{1}{m}{n-N\choose m}{n-N \choose m+1}$. As $N$ is a constant,
we have
\begin{equation} \lim_{n\rightarrow \infty} \sqrt[n]{
\frac{1}{m}{n-N\choose m}{n-N \choose m+1}} \end{equation}

Now restrict our attention to the particular case when $m=\lfloor n/3 
\rfloor $. We claim that permutations of this particular type are sufficiently
numerous to provide the lower bound we need.
 Using Stirling's formula, a routine computation yields that
in this case, we have
  \[\lim_{n\rightarrow \infty} \sqrt[n]{{n\choose m}{n\choose m}}
=\left(\frac{3^n}{2^{2n/3}}\right )^2 \geq 1.88^{2n}.
\]
Besides, we have more than $7.99^{2n/3}$ choices for the block-structured
string $S$ by which we replace $S^{dec}$. Therefore, we have proved the
following lower bound.

\begin{lemma} \label{lower}
For $n$ sufficiently large, the number of $n$-permutations of length $n$ 
that avoid the pattern 12453 is larger than
\[1.88^{2n} \cdot 7.99^{2n/3} \geq 14.12^n.\]
\end{lemma}

Lemma \ref{lower} and Corollary \ref{upper} together immediately yield
the following.

\begin{theorem} We have \[14.12 \leq L(12453)
\leq 14.66.\]
In particular, $L(12453)=\lim_{n\rightarrow \infty}
 \sqrt[n]{S_n(12453)}$ is not an integer.
\end{theorem}

\section{The exact value of $L(12453)$}

If we are a little bit more careful with our choice of $m$ in the argument
of the previous section, we can find the exact value of  
$L(12453)$. It turns out to be the upper bound
proved in Corollary \ref{upper}. 

\begin{theorem} \label{main} We have $L(12453)=(1+\sqrt{8})^2=9+4\sqrt{2}$.
\end{theorem}

\begin{proof} The above argument works for
any $m<n$ instead of $m=\lfloor n/3 \rfloor $,
 and for any positive real number $8-\epsilon
<8$ instead of 7.99. 

In order to find the best possible lower bound for $L(12453)$, set 
 $m=\alpha n$. Repeating the argument of the proof of Lemma \ref{lower},
we see that for $n$ sufficiently large, we have
 $S_n(12453)\geq {n\choose \alpha n}^2 (8-\epsilon) ^{n-\alpha n}. $
Then  the function
$f(\alpha)=\lim_{n\rightarrow \infty} 
({n\choose \alpha n}^2 (8-\epsilon) ^{n-\alpha n})^{1/n}$
 has a maximum
on the compact interval $[0,1]$. Choose the $\alpha $ providing that 
maximum. We claim that for
 that optimal $\alpha$, we must have $f(\alpha)=9+4\sqrt{2}$.
Indeed, we have
\[(1+\sqrt{8-\epsilon})^{2n} = \left 
(\sum_{m=0}^n {n\choose m}^2(8-\epsilon)^{n-m} \right )^2
 \leq (n+1)^2 
{n\choose \alpha n}^2 (8-\epsilon) ^{n-\alpha n},\]
because the square of an $(n+1)$-term sum has $(n+1)^2$ terms.
Taking $n$th roots, and then taking limits as $n$ goes to infinity, we
see that \[(1+\sqrt{8-\epsilon})^2\leq f(\alpha)\] for any
positive $\epsilon$, proving our claim. 
\end{proof}

\section{Some generalizations}

In this Section, we will provide some interesting
generalizations of our results. We will need the following simple recursive
properties of pattern avoiding permutations. 

\begin{proposition} \label{recprop}
Let $q$ be a pattern  of length $k$ that starts with 1,
 and let $q'$ be the
pattern of length $k+1$ that is obtained from $q$ by adding 1 to each entry
of $q$ and prepending it with 1. Let $p$ be a permutation whose string
of remaining entries is $S$. Then the following hold.
\begin{enumerate} 
\item If $S$ avoids $q$, then $p$ avoids $q'$.
\item If  $q$ itself starts with 1, then $p$ avoids $q'$ if and only if
$S$ avoids $q$. 
\end{enumerate} 
\end{proposition}

Iteratively applying part 2 of Proposition \ref{recprop}, and
 the method explained in the previous sections, we get the
following theorem.
\begin{theorem} \label{general} Let $k\geq 4$, and let
 $q_k$ be the pattern $12\cdots (k-3)(k-1)k(k-2)$. So $q_4=1342$, 
$q_5=12453$, and so on. Then we have
\[L(q_k)=(k-4+\sqrt{8})^2.\] 
 \end{theorem}

\begin{proof} Induction on $k$. For $k=4$, the result is proved in
\cite{bona}, and for $k=5$, we have just proved it in the previous section.
Assuming that the statement is true for $k$, we can prove the statement
for $k+1$ the very same way we proved it for $k=5$, using the result for
$k=4$, and part 2 of Proposition \ref{recprop}.
\end{proof}

The method we used to prove Lemma \ref{lupper} can also be used to prove the
following recursive result.

\begin{lemma} \label{gener}
Let $q$ be a pattern of length $k$ that starts with 1,
 and let $q'$ be the
pattern of length $k+1$ that is obtained from $q$ by adding 1 to each entry
of $q$ and prepending it with 1. Let $c$ be a
constant so that $S_n(q)<c^n$ for all $n$. Then we have
\[S_n(q')<(1+\sqrt{c})^{2n}=(1+c+2\sqrt{c})^n.\]
\end{lemma}

This is an improvement of the previous best result \cite{bonarec}, that
only showed $S_n(q')<(4c)^n$. 

The following generalization of Theorem \ref{main} can be proved just as
that Theorem is. 

\begin{theorem} \label{inter}
Let $q$ and $q'$ be as in Lemma \ref{gener}.  Then we have 
\[ L(q') =1+L(q)+2\sqrt{L(q)}.\]
\end{theorem} 

In a sense, this result generalizes Regev's result \cite{regev} that
showed that $L(12\cdots k)=(k-1)^2$. Our result shows that this particular
growth rate, that is, that $\sqrt{L(q)}$ grows by one as the pattern grows
by one, is not limited to monotone patterns. 

An interesting consequence of this Theorem is that if $q$ is as above, and
$L(q)<(k-1)^2$, in other words, $q$ is harder (or easier, for that matter)
 to avoid than the monotonic
pattern of the same length, then repeatedly prepending $q$ with 1 will not
change this. That is, the obtained new patterns will still be more difficult
to avoid than the monotonic pattern of the same length. 

Are the methods presented in this paper useful at all if the pattern $q$ does
not start in the entry 1? We will show that for most patterns $q$, the 
answer is in the affirmative, as far as a {\em lower} bound is concerned.
Let us say that the pattern $q$ is {\em indecomposable} if it cannot be 
cut into two parts so that all entries on the left of the cut are larger than
all entries on the right of the cut. For instance, 1423 and 3142 are 
indecomposable, but 3412 is not as we could cut it after two entries. 
Therefore, we call 3412 {\em decomposable}. It is
routine to verify that as $k$ grows, the ratio of indecomposable patterns
among all $k!$ patterns of length $k$ goes to 1. 

\begin{theorem} \label{indec} 
Let $q$ be an {\em indecomposable} pattern of length
$k$, and let  $L=\lim_{n\rightarrow 
\infty} \sqrt[n]{S_n(q)}$. Let $q'$ be defined as in Lemma \ref{gener}. 
Then we have \[ \lim_{n\rightarrow 
\infty} \sqrt[n]{S_n(q')} \geq 1+L+2\sqrt{L}.\]
\end{theorem}

\begin{proof}
This Theorem can be proved as Lemma \ref{lower}, and Theorem \ref{main}
are. Indeed, as $q$ is
indecomposable, any block-structured string $S$ will avoid $q$ if each
block does. Now apply part 1 of Proposition \ref{recprop} to see that
 our argument will still provide the required lower bound.  
\end{proof}

If $q=q_1q_2\cdots q_k$ is an indecomposable pattern, then so is its reverse
complement, that is, the pattern $q^{rc}$ whose $i$th entry is $k+1-q_{k+1-i}$
for all $i$. This leads to the following Corollary.

\begin{corollary}
Let $q$ be an indecomposable pattern, and let $q''$ be the pattern obtained
from $q$ by prepending $q$ with a 1, and appending the entry $k+2$ to the
end of $q$. Then we have
\[L(q'')\geq \left (2+\sqrt{L(q)} \right )^2.\]
\end{corollary}

Our methods will not provide an upper bound for $\lim_{n\rightarrow 
\infty} \sqrt[n]{S_n(q')}$ because the string $S$ of remaining entries
of a $q'$-avoiding permutation
does not have to be $q$-avoiding. (Only part 1 and not part 2 of Proposition
\ref{recprop} applies.) That condition is simply
sufficient, but not necessary, in this general case. Nevertheless, Theorem
\ref{indec} is interesting. It shows that for almost all patterns $q$, if
we prepend $q$ by the entry 1, the limit of the corresponding Stanley-Wilf
sequence will grow at least as fast as for monotone $q$. If $q$ started in 
1, then this growth will be the same as for monotone $q$. 

Now it is a little easier to understand why, in the case of length 4, the
patterns that are the hardest to avoid, are along with certain equivalent
ones, 1423 and 1342. Indeed, removing the starting 1 from them, we get
the {\em decomposable} patterns 423 and 342. As these patterns are
 decomposable, Theorem \ref{indec} does not hold for them, so the limit
of the Stanley-Wilf sequence for the patterns 1423 or 1342 does not have
to be at least $1+4+4=9$, and in fact it is not.

A particularly interesting application of Theorem \ref{inter} is as follows.
Recall that a {\em layered} pattern is a pattern that consists of
decreasing subsequences (the layers) so that the entries increase
among the layers. For instance, 3217654 is a layered pattern.
In 1997, several people (including present author) 
have observed, using numerical evidence
computed in \cite{West}, that if $q$ is a layered pattern of length $k$,
 then for small $n$, the
inequality $S_n(12\cdots k)\leq S_n(q)$ seems to hold.
 We will now show that this
is indeed true in the sense of logarithmic asymptotics.  

\begin{theorem} \label{layered}
Let $q$ be a layered pattern of length $k$. Then we have
\[L(q)\geq (k-1)^2.\]
Equivalently, $L(q)\geq L(12\cdots k)$.
\end{theorem}

In order to prove Theorem \ref{layered}, we need the following powerful
Lemma,
due to Backelin, West, and Xin.

\begin{lemma}  \label{bwx} \cite{bwx}
Let $v$ be any pattern of length $k-r$. Then for all positive integers $n$ and
$r<k$,
we have \[S_n(12\cdots r v)=S_n(r(r-1) \cdots 21v),\]
where $v$ is taken on the set $\{r+1,r+2,\cdots ,k\}$. 
\end{lemma}

Now we are in position to prove Theorem \ref{layered}.

\begin{proof} (of Theorem \ref{layered}.)
Induction on $k$.  If $q$ has only one layer, then
$q$ is the decreasing pattern, and the statement is obvious. Now assume 
$q$ has at least two layers, and that we know the statement for all layered
patterns of length $k-1$. As $q$ is
layered, it is of the form $r(r-1)\cdots 21v$ for some $r$, and some
layered pattern $v$. Therefore, Lemma \ref{bwx} applies, and we have
$S_n(q)=S_n(12\cdots rv)$. If this last pattern is denoted
by $q^*$, then we obviously also have $L(q)=L(q^*)$. We further
denote by $q^{*-}$ the pattern obtained from $q^{*}$ by removing its
first entry. Note that $q^{*-}$ is still a layered pattern, just its first
several layers may have length 1. 

Assume first that $r>1$. Then note that $q^{*-}$ starts with its smallest
entry. Therefore, 
Theorem \ref{inter} applies,  and by the  induction hypothesis we have
\[L(q)=L(q^*)=1+L(q^{*-})+2\sqrt{L(q^{*-})}\geq 1 + (k-2)^2+2(k-2)=
(k-1)^2,\]
which was to be proved.

Now assume that $r=1$. Then  $q$ is a layered
pattern that starts with a layer of length 1. Therefore, instead of
applying Theorem \ref{inter}, we need to, and almost always can, 
apply Theorem \ref{indec} for the pattern $q^{*-}$. Indeed, 
$q^{*-}$ is a layered pattern, and as such, is indecomposable, except when
it has only one layer, that is, it is the decreasing permutation.
Therefore, Theorem \ref{indec} implies
\[L(q)\geq 1+L(q^{*-})+2\sqrt{L(q^{*-})}\geq 1 + (k-2)^2+2(k-2)=
(k-1)^2.\]

Finally, if $r=1$, and $q^{*-}$ is the decreasing pattern, then we simply have
$q=1k(k-1)\cdots 2$. In that case, our statement is just a special case
of Lemma \ref{bwx}. Indeed, choosing $v$ to be the decreasing pattern, 
Lemma \ref{bwx} shows $S_n(q)=S_n(12\cdots k)$.

This completes the proof.
\end{proof}

Here is another way in which  our results start explaining why certain 
patterns are easier to avoid than others. We formulate our observations in
the following Corollary. 

\begin{corollary} 
Let $q_1$ and $q_2$ be patterns so that 
$L(q_1)\leq L(q_2)$. Let 
$q_i'$  be the pattern obtained from $q_i$ by prepending $q_i$ by a 1. 
Furthermore, let $q_1$ start with the entry 1,
 and let $q_2$ be indecomposable.
Then we have
\[L(q_1') =1+L(q_1)+2\sqrt{L(q_1)} \leq
 1+L(q_2)+2\sqrt{L(q_2)} \leq
L(q_2') .\]
\end{corollary}  

For instance, if we set $q_1=123$ and $q_2=213$, we get the well-known
statement weakly comparing the limits of the Stanley-Wilf sequences
of 1234 and 1324, first proved in \cite{bonarec}.

\section{Further Directions}
Our results raise two interesting kinds of questions. We have seen that the
limit of a Stanley-Wilf sequence is not simply not always an integer, but 
also not always rational. Is it always an {\em algebraic} number? If yes,
can its degree be arbitrarily high? Can it be more than two? Is it always
an {\em algebraic integer}, that is, the root of a {\em monic} 
polynomial with integer
coefficients? The results so far leave that possibility open.

The second question is related to the size of the limit of $\sqrt[n]{S_n(q)}$
if $q$ is of length $k$. The largest value that this limit is known to take
is $(k-1)^2$, attained by the monotonic pattern. Before present paper, 
the smallest known value, in terms of $k$, for this limit was $(k-1)^2-1=8$,
attained by $q=1342$. As Theorem \ref{general} shows, the value
$(k-4+\sqrt{8})^2$ is also possible. As $k$ goes to infinity, the
{\em difference} of  the assumed maximum $(k-1)^2$ 
and this value also goes to infinity, while their ratio goes to 1. 
Is it possible to find a series of patterns $q_k$  so that this ratio does
not converge to 1?  We point out that it follows from a result of 
P. Valtr (published in \cite{kaiser})
 that for any pattern $q$ of length $k$, we have 
$\lim_{n\rightarrow \infty} \sqrt[n]{S_n(q)} \geq e^{-3}k^2$, so the
mentioned ratio cannot be more then $e^3$.

Finally, now that the Stanley-Wilf conjecture has been proved, and we know
 that the limit of a Stanley-Wilf sequence always exists, we can ask what
the largest possible value of this limit is, in terms of $k$. In 
\cite{arratia}, R. Arratia conjectured that this limit is at most $(k-1)^2$,
and, following the footsteps of Erd\H{o}s, he offered 100 dollars for a
proof or disproof
 of the conjecture $S_n(q)\leq (k-1)^{2n}$, for all $n$ and $q$.
Our results provide some additional support for this conjecture as they 
show that there is a wide array of patterns $q$ for which $\sqrt{L(q)}$ grows
by one when $q$ is prepended by the entry 1. In fact,
numerical evidence suggests that even the following stronger version of
Arratia's conjecture could be true. 

\begin{conjecture} Let $q$ be a pattern of length $k$. Then $L(q)\leq 
(k-1)^2$, where equality holds if and only if $q$ is layered, or the
reverse of $q$ is layered.
\end{conjecture}

\centerline{\bf Acknowledgment} 

I am indebted to the anonymous referee whose careful reading significantly
improved the presentation of my results.

\end{document}